\documentclass{article}

\newtheorem{thm}{Theorem}
\newtheorem{lem}{Lemma}

\newtheorem{prop}{Proposition}

\newcommand{\pa}{\partial}

\newcommand{\idsum}{\sum_{i=1}^d}
\newcommand{\jdsum}{\sum_{j=1}^d}
\newcommand{\kdsum}{\sum_{k=1}^d}
\newcommand{\ijdsum}{\sum_{i, j=1}^d}
\newcommand{\adsum}{\sum_{\alpha = 1}^d}
\newcommand{\deli}{\frac{\partial}{\partial x^i}}
\newcommand{\delj}{\frac{\partial}{\partial x^j}}
\newcommand{\dv}{{\rm div}}
\newcommand{\intt}{\int_0^T}
\newcommand{\tilU}{\tilde{U}}

\newcommand{\tilB}{\tilde{B}}
\newcommand{\tiltau}{\tilde{\tau}}
\newcommand{\tr}{{\rm tr}}

\begin{document}
\title{Stochastic analysis in a tubular neighborhood\\
{\Large or Onsager-Machlup functions revisited}\\
{\large (On the occasion of 70th birthday of Y.Takahashi)}}
\author{Keisuke Hara\thanks{Email: hara.keisuke@gmail.com} \and Yoichiro Takahashi}
\date{}
\maketitle

\section{A personal introduction by K.H.}

Inspired by K.~It\^o's work on Feynman path integrals
and discussions on the effect of curvature among physicists,
Y. Takahashi started studying the problem of the most probable path
(or the Onsager-Machlup function) in late 1970's
and completely determined the function
in his joint work with S.~Watanabe \cite{TW} (1980).

Soon after their work, T.~Fujita and S.~Kotani \cite{FK} (1982)
gave the general result.
While they used a singular perturbation method for P.D.E.,
Takahashi-Watanabe's proof is almost probabilistic
except the only one analytic part,
which treats asymptotic behaviour of a certain Wiener functional.
Though this part gives a key observation,
Takahashi removed it in a joint work with me
to get a purely probabilistic proof and to extend the result
to pinned diffusion processes \cite{HT} (1996).

On the other hand, he had been studying the Wiener functional
in his way, and wrote {\it a draft}.
In Takahashi-Watanabe \cite{TW},
we can see how the scalar curvature
term appears through the functional,
while we cannot see the effect in \cite{HT}
because the proof is too simple.
Actually it appears from a kind of ergodic theorem
due to the rapid rotation of the spherical motion
as the tube is shrinking.
His project was to study the effect
with pure stochastic analysis.

Unfortunately, the draft had been missing for long years because
he gave up the work to move to the other projects and
also because ``yet another" proof became somehow needless
after the third proof.

However, recently I found the draft (written by himself
with a primitive word processing software) in my library.
Though it is hard to read because of
his confusing writing, missing many parts, and the oxidization,
I also found it unexpectedly interesting and tricky.
Then, I tried to restore it as much as possible, reassembled,
and I submit it as our {\it joint work}
on the occasion of his 70th birthday
after 20 years since he handed me the draft.

\section{Introduction and the main result}

Let $x(t)$ be a non-degenerate smooth diffusion process
on a $d$-dimensional Riemannian manifold $M$
and consider the asymptotics of the probability of
the trajectory being confined in a small tubular
neighbourhood of a prescribed smooth curve $\gamma$ on $M$.
The second author and S.~Watanabe \cite{TW} proved
that the following asymptotics is true.
The functional $L$ is the Onsager-Machlup function,
which is similar to the Lagrangian. 
The point is the effect of the curvature
to the ``most probable path".

\begin{thm}\label{thm:Main}
\[
\lim_{\delta \to 0}
\frac{P_{\gamma(0)} \left[ \max_{0 \leq t \leq T} d(x(t), \gamma(t)) \leq \delta \right]}
 {P_0 \left[ \max_{0 \leq t \leq T} |B(t)| \leq \delta \right]}
= \exp \{ -S(\gamma) \},
\]
where
$B(t)$ is the standard Brownian motion in the $d$-dimensional Euclidean space,
$d( \cdot, \cdot)$ is the Riemannian distance on $M$, 
\[
S(\gamma) = S_T(\gamma) = \int_0^T L(\gamma(t), \dot{\gamma}(t) ) dt,
\]
\[
L(x, v) = \frac{1}{2} |f - v|_x^2 + \frac{1}{2} \dv f(x)
 + \frac{1}{12} R(x),
\]
and $R(x)$ is the Scalar curvature at $x$.
\end{thm}

The theorem above was proved firstly by Takahashi-Watanabe \cite{TW}.
In their almost probabilistic proof,
they studied asymptotic behaviour of a certain Wiener functional
to produce a key geometric quantity as an ergodic effect.
Though this analysis is crucial, it is the only one
``non-probabilistic" part in their proof.
Soon after their work, Fujita-Kotani \cite{FK} gave
an analytic proof to the theorem
by a singular perturbation method for P.D.E.

Following these two proofs, the authors \cite{HT} gave another proof,
which is simple and purely probabilistic, 
and extended the result to pinned diffusion processes.
However, in return for the simplicity, the new proof lost
a view to the ergodic effect.
The aim of this paper is to give yet another (purely) probabilistic proof
by recovering such precise study on the Wiener functional.

\bigskip

To highlight the difference between the old proofs and our new proof,
here we recall the argument used in \cite{TW}
and sketch our new idea.

The first key point is the existence of the Lagrangian.
It follows from the two facts.
The first fact is that for any given curve $\gamma$
the law of the distance process $d(x(t), \gamma(t))$
is absolutely continuous to the law
of a Bessel process with index $d = {\rm dim} M$
if they are restricted to a sufficiently small neighbourhood
of the curve.
We can find the law by adjusting the drift term
expressed in normal coordinates along the curve $\gamma$,
which turns out to be the difference of the Coriolis drift
and the Besselization drift.
Thus we find a diffusion process whose radial part admits the Bessel process
and take it as the reference process to compute
the Radon-Nykodim density given by Girsanov formula.

The latter fact is that the spherical motion of the reference
process is governed by Levy's stochastic areas,
which are orthogonal to the martingales
governing the radial motion.
Hence we can apply Kunita-Watanabe theorem
to factorize the density into the Lagrangian
and the asymptotically vanishing term.

The second key point is the derivation of the form of the
Lagrangian $L(x, v)$, which consists of three terms as
its stated above.
The first term  $ \frac{1}{2} |v-f|_x^2$
comes out from the quadratic variation term in the Girsanov exponent.
The other terms are obtained as the correction
due to the transformation of the martingale term
in the exponent into the stochastic linear integral.
We apply the stochastic Stokes theorem to estimate the line integral.

The third term is the mean curvature term, which is hardest to obtain.
Though it appears similarly as the difference of the quantities,
it can be obtained more indirectly as an ergodic phenomenon.
In Takahashi-Watanabe \cite{TW}, this part (and only this part)
was proved by appealing to an analytical method,
i.e., the singular perturbation technique,
using which Fujita-Kotani \cite{FK} gave another whole proof.
After these two proofs, the authors gave a purely probabilistic method
by ``smooth Besselization technique" \cite{HT},
which was introduced  in \cite{H} by the first author.
However, it lost the concrete study of the ergodic effect
in return for simplicity.

In this paper we prove this part by stochastic analysis.
Precisely to say, the word ``ergodic" stated above
is not one in the usual sense but it is asymptotic.
We can find a process that may
be called asymptotically tangent to the original process
as the tube shrinks,
which is nothing but a Brownian motion.
To justify this argument,
we need to prepare a {\it strong} comparison theorem
between diffusion processes with different diffusion constants.
The comparison theorem guarantees that
an asymptotic ergodic theorem is inherited by the original diffusion
process from the asymptotically tangent process
(i.e., a Brownian motion).

The trick for the strong comparison is to adjust the Brownian motions in
the stochastic differential equations
such that the distance between the two processes is kept small enough.
We prepare a key lemma of linear algebra for this task.

\section{Coriolis drift and Besselization drift}

First we prepare necessary geometric quantities and
introduce two important vectors: Coriolis drift and Besselization drift.
Though we use properties of normal coordinates
and basic asymptotics of the geometric tensors in the coordinates,
we omit the proofs.
They can be found in textbooks like \cite{W} or \cite{G}.

The equation associated with the diffusion process $x(t)$ is,
by the assumption, a non-degenerate second order
parabolic differential equation on $M$ as follows:
\[
\frac{\partial u}{\partial t} = \frac{1}{2} \triangle_g u + fu,
\]
where
$\triangle_g$ is the Laplace-Beltrami operator
for Riemannian metric $g$
and
$f$ is a vector field on $M$.
In local coordinates they are expressed as
\[
\triangle_g u(x) = \frac{1}{\sqrt{g(x)}} \ijdsum
\frac{\pa}{\pa x^j}
 \left(
    \sqrt{g(x)} g^{ij}(x) \frac{\pa u}{\pa x^j}(x)
 \right),
\quad
f u(x) = \idsum f^i (x) \frac{\pa u}{\pa x^i} (x).
\]
Here  $(g^{ij} (x))$ is the diffusion coefficient
in the coordinate, and its inverse matrix $(g_{ij} (x))$
gives the metric
$ds^2 = \sum g_{ij} (x) dx^i dx^j$.
We also used the conventional notation $g(x) = \det(g_{ij}(x))$.
The norm $| \cdot |_x$ is the Riemannian norm on the tangent
space $T_x M$ at the point $x$, i.e.,
\[
|v|_x^2 = \ijdsum g_{ij} (x) v^i v^j
\quad \mbox{for} \quad v = \idsum v^i \frac{\pa}{\pa x^i}.
\]
The divergence of a vector field $f = \sum f^i (\partial / \partial x^i)$
is denoted by $\dv f$, i.e.,
\[
\dv f(x) = \frac{1}{\sqrt{g(x)}} \idsum \deli \left(
    \sqrt{g(x)} f^i(x) \right).
\]

We rewrite the Laplace-Beltrami operator as follows:
\[
\frac{1}{2} \triangle_g u(x)
= \ijdsum \frac{1}{2} g^{ij}(x) \deli \delj u(x)
 + \idsum a^i (x) \deli u(x)
\]
where $a(x) = (a^i (x))$ is defined by
\[
a^i(x) = \jdsum \frac{1}{\sqrt{g(x)}} \delj \left(
    \sqrt{g(x)} g^{ij} (x) \right).
\]
This drift term expressed by $a(x)$ is called Coriolis drift.

Let $\gamma(t) \, (0 \leq t \leq T)$
be a smooth curve on the manifold $M$.
In the following, 
we take the lifting along the curve $\gamma$, i.e.,
normal coordinates along $\gamma$ or Fermi coordinates (See \cite{G}).

To write down our stochastic differential equation (S.D.E.),
we prepare the square root $\sigma (t,x)$ of the matrix
$(g^{ij} (t, x))$:
\[
\sigma(t, x)
= (\sigma^{ij} (t, x))_{i,j = 1, \dots, d}
= (\sigma_\alpha^i (t, x))_{i, \alpha = 1, \dots, d}
= (g^{ij} (t, x))^{1/2}.
\]

Since $\sum_j g^{ij}(x) x^{j} = x^i$ holds in normal coordinates,
we have $\sum_\alpha \sigma_\alpha^i (x) x^{\alpha} = x^i$.
Note that we have also the same identities in Fermi coordinates,
i.e.,
$\sum_j g^{ij}(t, x) x^j = x^i$ and
$\sum_\alpha \sigma_\alpha^i(t, x) x^{\alpha} = x^i$.

\begin{lem}
There exists a local vector field $c(t,x)$ defined
in normal coordinates along $\gamma$ such that
the solution $Y(t)$ to the S.D.E.
\begin{equation}\label{eqY}
dY(t) = \sigma(t, Y(t)) dB(t) + c(t, Y(t)) dt
\end{equation}
admits the Bessel process of index $d$ as its radial motion.
Moreover, if, in addition, one requires the symmetric condition
$ c^i (t, x) x^j = c^j (t, x) x^i \, (i, j = 1, \dots, d)$,
then it is uniquely determined and is given as
\[
c^i (t, x) = \frac{x^i}{2 |x|^2} \jdsum (1 - g^{jj} (t,x)).
\]
\end{lem}

{\it Remark.}
The choice of this Besselization drift is a key point.
The drift $c(t,x)$ has the coordinate symmetry,
but it has also the singularity at $x = 0$.
On the other hand, in \cite{HT} (originated in \cite{H})
we broke the symmetry to choose another Besselization drift,
which is totally smooth without any singularity.

{\it Proof.}
Apply It\^o formula to  $\jdsum Y^j (t)^2$
and compare it with S.D.E. of the Bessel process.
See \cite{TW} for details.

\bigskip

We have the asymptotic expansion of our geometric tensors
$\sigma(t, x)$, $a(t, x)$, and $c(t, x)$ as $x \to 0$ in
normal coordinates along the curve as follows.

\begin{lem} \label{lem:del-a}\label{lem:del-c}
\begin{eqnarray*}
\sigma^{ij} (t, x)
 &=& \delta^{ij} + \frac{1}{6} R_{ikjl}(t, 0) x^k x^l + O(|x|^3),\\
\jdsum \frac{\pa}{\pa x^j} a^j (t, x)
 &=& -\frac{1}{3} R(t, 0) + O(|x|^2),\\
\jdsum \frac{\pa}{\pa x^j}  c^j (t, x)
 &=& \frac{d}{6} \ijdsum R_{ij}(t, 0) u^i u^j + O(|x|),
\end{eqnarray*}
where $\delta^{ij}$ is the Kronecker delta,
$R_{ijkl}$ is the Riemannian curvature,
$R_{ij}$ is the Ricci curvature,
$R$ is the scalar curvature,
and $u^i = x^i / |x^i|$.
\end{lem}

\section{Girsanov formula}

In the same way as the preceding papers \cite{TW} \cite{FK} \cite{HT}, 
we write the S.D.E. of the process $X(t)$ obtained
by lifting the original process $x(t)$
as follows.

\begin{equation}\label{eqX}
dX^i (t) = \adsum \sigma_\alpha^i (t, X(t)) dB^\alpha (t)
 + \tilde{b}^i (t, X(t)) dt,
\end{equation}
where the drift term $\tilde{b}$ is
\[
\tilde{b}^i(t, x)
= a^i(t, x) + f^i(t, x) - \dot{\gamma}^i(t)
\]
in normal coordinates along $\gamma(t)$.
To use it later, we set
\[
b^i (t, x) = f^i(t, x) - \dot{\gamma}^i(t).
\]

Let us compare $X(t)$ with  $Y(t)$
whose radial part is a Bessel process.
By Girsanov formula, we obtain the following
(For the proof, see \cite{TW} or \cite{HT}).

\begin{prop}\label{prop:Girsanov}
As $\delta \to 0$,
\begin{eqnarray*}
&& P_{\gamma (0)} \left[
\max_{0<t<T} d(x(t), \gamma(t)) < \delta \right] 
= P_0 \left[ \max_{0<t<T} d(X(t), 0) < \delta \right]\\
&=& \exp \{ -S_T (\gamma) + O(\delta) \}
E_0 \left[ \exp \left( \int_{Y[0, T]} \alpha + \int_0^T \beta (t) dt
 \right); |Y|_T < \delta \right],
\end{eqnarray*}
where $|Y|_T = \max_{[0, T]} |Y(t)|$,
the integral over $Y[0, T]$ is the stochastic line integral
over the space-time curve $(t, \gamma(t)), 0 \leq t \leq T$,
the integrand $\alpha$ is the space-time 1-form
(which is degenerated in $t$-direction)
given by
\[
\alpha = \idsum \alpha_i dx^i
\quad\mbox{with} \quad
\alpha_i = \jdsum (a^j (t,x) + b^j(t,x) - c^j(t,x)) g_{ij}(t,x),
\]
the integrand $\beta (t)$ stands for the function
\[
\beta (t) = \frac{d}{12} \ijdsum R_{ij} (\gamma (t))
\left( \frac{Y^i (t)}{|Y(t)|} \frac{Y^j (t)}{|Y(t)|}
 - \frac{\delta^{ij}}{d} \right).
\]
\end{prop}

\section{Stochastic Stokes theorem}

Now let us apply the stochastic version of the Stokes theorem
(\cite{IW}) to get the estimate of the stochastic line
integral of the 1-form $\alpha$ over the space-time
curve $Y[0,t]$.
Let $\Sigma$ be the (random) surface
$\{ sY(t) \}_{0 \leq s \leq 1, 0 \leq t \leq T}$.
Then its boundary $\pa \Sigma$ consists of four curves:
the trajectory $Y[0, T]$ (corresponding to $s = 1$),
the segment from $(T, Y(T))$ to $(T,0)$ (corr. to $t = T$),
the segment from $(T, 0)$ to $(0, 0)$ (corr. to $t = 0$),
and the degenerated segment consisting of one point $(0, 0)$
(corr. to $s = 0$).

Note that the line integrals over those segments are of
order $O(\delta)$ on the set where $|Y|_T < \delta$
because $d(sY(t)) = Y(t) ds$ on the segments
where $t$ is a constant.
Then, stochastic Stokes theorem says that
\[
\int_\Sigma d \alpha = \int_{\pa \Sigma} \alpha
= \int_{Y[0,T]} \alpha + O(\delta T).
\]

Therefore, the estimate of the line integral
is deduced to the study of the area,
which is expressed concretely as follows.

\begin{lem}
Denote the stochastic areas (L\'evy's area) by
\[
A^{ij} (t) = \int_0^t (Y^i \circ dY^j - Y^j \circ dY^i)
\]
and set
\[
\alpha_{ij} (t, x) = \frac{1}{4} \int_0^1 \left\{
 \pa_i \alpha_j (t, sx) - \pa_j \alpha_i (t, sx) \right\} ds.
\]
Then,
\[
\int_\Sigma d\alpha = \ijdsum \intt \alpha_{ij} (t, Y(t))
 \circ dA^{ij} (t).
\]
\end{lem}

{\it Proof.}
Obvious from the definitions.
(Notice that we omitted the factor $(1/2)$ from $A^{ij}$,
which is set into $\alpha_{ij}$.)

\bigskip

The following lemma is crucial,
because we need the orthogonality
to factorize the density into the Lagrangian term
and the asymptotically vanishing term by Kunita-Watanabe theorem.

\begin{lem}
The stochastic areas $A^{ij} (t)$ are martingales which are
orthogonal to the radial motion $|Y(t)|$.
In other words, the quadratic variation process vanishes:
\[
\langle A^{ij}, |Y| \rangle (t) = 0.
\]
\end{lem}

{\it Proof.}
Show 
 $(Y^i \circ dY^j - Y^j \circ dY^i) \kdsum Y^k dY^k = 0$
by stochastic calculus.
See \cite{TW} or \cite{HT}.

\bigskip

Therefore we proceed to the next step as follows:

\begin{prop}\label{prop:preMain}
For $p > 0$, set
\[
M_p(t) = p \int_0^t \alpha_{ij} (s, Y(s)) dA^{ij} (s)
 - \frac{p^2}{2} \int_0^t \alpha_{ij} \alpha_{kl}
 d \langle A^{ij}, A^{kl} \rangle (s),
\quad M(t) = M_1(t),
\]
and
\[
L(t) = \frac{d}{12} \int_0^t  R_{ij} (\gamma (s))
 \left\{ U^i (s) U^j (s) - \frac{\delta^{ij}}{d} \right\} ds,
\]
where $U(t)$ is the spherical part of $Y(t)$, i.e.,
$ U^i (t) = Y^i (t) / |Y(t)|$.
Then,
\begin{eqnarray*}
&&P \left[ \max_{0 < t < T} d (x(t), \gamma(t)) < \delta \right]\\
&=&
\exp \left\{ -S_T(\gamma) + O(\delta T) \right\} \,
 E \left[\exp \{ M(t) + L(t) \} \, ; \, |Y|_T < \delta \right].
\end{eqnarray*}
\end{prop}

{\it Proof.}
Note that on the set $|Y|_T < \delta$,
\[
\intt \alpha_{ij} \alpha_{kl} d \langle A^{ij}, A^{kl} \rangle (t)
 = O(\delta^2 T),
\]
because
$(Y^i \circ Y^j - Y^j \circ Y^i) (Y^k \circ Y^l - Y^l \circ Y^k) 
 = O(\delta^2) dt$.
And also,
\[
\langle \alpha_{ij}, A^{ij} \rangle (T) = O(\delta T),
\]
because $d \alpha_{ij}(Y^i \circ Y^j - Y^j \circ Y^i) = O(\delta) dt$.
Hence,
\[
\intt \alpha_{ij} \circ dA^{ij} = M(T) + O(\delta T)
\quad \mbox{on $|Y|_T < \delta$}.
\]
Consequently the assertion follows from Proposition \ref{prop:Girsanov}.

\bigskip

{\it Remark.}
If $(M, g)$ is an Einstein space,
we have the relation $R_{ij} = (1/d) g_{ij}$ and so $L(t) = 0$.
Therefore, applying Kunita-Watanabe theorem
to the exponential martingale part,
we already have proved Theorem \ref{thm:Main} in this case.

\section{Pathwise adjustment}

The following sections are the new part
where we study the Wiener functional $L(t)$ with stochastic analysis.

Let us recall that the process $Y(t)$ is governed by the S.D.E.:
\[
dY(t) = \sigma(t, Y(t)) dB(t)
 + \frac{Y(t)}{2|Y(t)|^2} \left(
  d - \jdsum g^{jj} (t, Y(t)) \right) dt.
\]

Since $\sigma(t, Y(t)) = I + O(|Y(t)|^2)$,
one may expect that the process $Y(t)$ converges as $\delta \to 0$
under $P[ \, \cdot \, | \, |Y(t)|_T < \delta]$
to a $d$-dimensional Brownian motion,
possibly with a drift,
in some sense and that we can estimate its convergence rate.
However,
usual comparison theorems for S.D.E.s
give very crude information when the diffusion coefficients
vary in contrast with comparison theorems for partial differential equations.
Moreover what we want to deal with is its spherical part
$U(t) = Y(t)/ |Y(t)|$.
We solve these problems as follows.

By It\^o formula, the spherical part $U(t)$
is governed by the S.D.E.:
\begin{eqnarray}
dU(t) &=& \frac{1}{|Y(t)|} \left( I - U(t) \otimes U(t) \right)
 \sigma (t, Y(t)) dB(t) \nonumber \\
 &&- \frac{1}{2|Y(t)|^2} \left( \jdsum g^{jj} (t, Y(t)) - 1 \right)
  U(t) dt, \label{eqU}
\end{eqnarray}
where we denote by $u \otimes v$
the matrix $(u^i v^j)_{1 \leq i, j \leq d}$.

On the other hand,
the spherical part $\tilde{U} (t)$ of a $d$-dimensional
standard Brownian motion $\tilde{B}(t)$ is,
of course, governed by the S.D.E.:
\begin{equation}\label{eqUtil}
d\tilde{U}(t) = \frac{1}{|\tilde{B}(t)|}
 \left( I - \tilde{U}(t) \otimes \tilde{U}(t) \right) d\tilde{B}(t)
 - \frac{d-1}{2 |\tilde{B}(t)|^2} \tilde{U}(t) dt.
\end{equation}

The idea to obtain a good comparison estimate between
(\ref{eqU}) and (\ref{eqUtil}) is to adjust the Brownian motions
so that the radial part $| \tilde{B} (t) |$ of $\tilde{B}(t)$
is {\it pathwise} equal to the radial part $|Y(t)|$ of $Y(t)$.
With this trick,
we can manage the fluctuation of the spherical part
enlarged by the fluctuation of radial parts.
For this sake we prepare the following lemma.

\begin{lem}\label{lem:linalg}
Let $n = 1 + d(d-1)/2$.
Then there are linear maps $J^i \, (i = 1, \dots, d)$
from ${\bf R}^d$ to  ${\bf R}^n$ with the following properties:

{\rm (a)}
Each $J^i$ is an isometry of ${\bf R}^d$ into  ${\bf R}^n$.
In other words, we have
\[
\langle J^i u, J^i v \rangle_n = \langle u, v \rangle_d
\]
for any $u, v$ in ${\bf R}^d$,
where $\langle \, \cdot \, , \, \cdot \, \rangle_n$
and $\langle \, \cdot \, , \, \cdot \, \rangle_d$
denote the inner products in ${\bf R}^n$ and ${\bf R}^d$ respectively.

{\rm (b)}
For each unit vector $u \in {\bf R}^d$,
$\{J^i u \} \, (1 \leq i \leq d)$ are orthonormal, i.e.,
\[
\langle J^i u, J^j u \rangle_n = \delta^{ij}.
\]

{\rm (c)}
There is a unit vector $e_0 \in {\bf R}^n$ such that
\[
\idsum u^i J^iu = |u|_d^2 e_0
\]
for every $u = (u^i)_{1 \leq i \leq d}$ in ${\bf R}^d$.
\end{lem}

{\it Proof.}
Let $e_\alpha \, (1 \leq \alpha \leq d)$ be the canonical basis
in ${\bf R}^d$ and $e_0, e_{\alpha \beta} \, (\alpha < \beta)$
be an orthonormal basis in ${\bf R}^n$.
(Thus, $n = 1 + d(d-1)/2$.)
Set
\[
J^i e_\alpha
= \left\{
 \begin{array}{ll}
 e_{i \alpha} & \mbox{if $i < \alpha$},\\
 e_0 & \mbox{if $i = \alpha$},\\
 -e_{\alpha i} & \mbox{if $i > \alpha$}.
 \end{array}
 \right.
\]
Then a direct, elementary computation shows that
\[
\langle J^i e_\alpha, J^j e_\beta \rangle
 = \delta_{ij} \delta_{\alpha \beta} - \delta_{i \beta} \delta_{\alpha j}
 + \delta_{\alpha i} \delta_{j \beta}.
\]
Now from the equation above  one obtains (a) by putting $i = j$ and
(b) by putting $\alpha = \beta$.
Finally,
\begin{eqnarray*}
\idsum u^i J^i u
&=& \idsum \adsum u^i u^\alpha J^i e_\alpha\\
&=& \sum_{i < \alpha} u^i u^\alpha e_{i \alpha}
 + \sum_i (u^i)^2 e_0
 + \sum_{i > \alpha} - u^i u^\alpha e_{\alpha i}
= |u|^2 e_0.
\end{eqnarray*}
Hence we obtain (c).

\bigskip

Now let us take a $(1 + d(d-1)/2)$-dimensional Brownian motion $W$
and consider the following S.D.E.s:
\begin{eqnarray*}
dY(t) &=& \sigma(t, Y(t)) \langle JU(t), dW(t) \rangle
 - \frac{1}{2|Y(t)|^2} \jdsum \left\{
  g^{jj} (t, Y(t)) - 1 \right\} Y(t) dt,\\
d\tilde{Y}(t) &=& \langle J \tilde{U}(t), dW(t) \rangle,
\end{eqnarray*}
where
$U(t) = Y(t)/|Y(t)|$, $\tilde{U}(t) = \tilde{Y}(t)/|\tilde{Y}(t)|$,
and $\langle Ju, v \rangle$ stands for the $d$-dimensional vector
whose components are $\langle J^i u, v \rangle_n \, (1\leq i \leq d)$.
Solve these S.D.E.s above
and then we can define the desired Brownian motions
$B(t)$ and $\tilde{B}(t)$ by
\[
dB(t) = \langle JU(t), dW(t) \rangle
\quad \mbox{and} \quad
d{\tilde{B}}(t) = \langle
 J\tilU(t), dW(t) \rangle.
\]
Then it is obvious that
\[
\langle U(t), dB(t) \rangle
= \langle \tilde{U}(t), d\tilde{B}(t) \rangle = dW_0,
\]
where $W_0(t)$ is the 1-dimensional Brownian motion
which is defined by
\[
W_0(t) = \int_0^t \langle e_0, dW(s) \rangle.
\]
With these Brownian motions we obtain the desired S.D.E.s
(\ref{eqU}) and (\ref{eqUtil}) for which,
by the property (c) of Lemma \ref{lem:linalg},
the radial parts coincide with each other if the initial data do so.
In fact,
\begin{eqnarray*}
d( |Y(t)|^2 ) &=& 2 \langle Y(t), dY(t) \rangle
 + \langle dY(t), dY(t) \rangle\\
&=& 2 \langle Y(t), \sigma (t, Y(t)) dB(t) \rangle\\
 && -2 \langle Y(t), \frac{1}{2|Y(t)|^2} \jdsum (g^{jj} - 1) Y(t) dt \rangle
 + \jdsum g^{jj} dt\\
&=& 2 \langle Y(t), dB(t) \rangle + d \cdot dt\\
&=& 2 |Y(t)| \langle U(t), dB(t) \rangle + d\cdot dt
= 2|Y(t)| dW_0 (t) + d \cdot dt.
\end{eqnarray*}
On the other hand,
\begin{eqnarray*}
d(|\tilde{Y}(t)|^2)
&=& 2 \langle \tilde{Y} (t), d \tilde{B} (t) \rangle
 + d \cdot dt
= 2 | \tilde{Y} (t) | \langle \tilde{U}(t), d\tilde{B}(t) \rangle
 + d \cdot dt\\
&=& 2 |\tilde{Y}(t) | dW_0 (t) + d \cdot dt.
\end{eqnarray*}
Hence, $|Y(t)| = |\tilde{Y}(t)|$ if $|Y(0)| = |\tilde{Y}(0)|$.

\section{The inner product of the spherical parts}

Next let us deduce the S.D.E. for the inner product
of $U(t)$ and $\tilU (t)$.

\begin{lem}\label{lem:uu}
The inner product $\langle U(t), \tilU(t) \rangle$
satisfies the following S.D.E.:
\[
d \langle U, \tilU \rangle
= R(t) H(t) dW_1(t) - \frac{1}{2} R(t)^2 G(t) \langle U, \tilU \rangle dt,
\]
where $W_1(t)$ is a 1-dimensional Brownian motion
that is independent of $W_0(t)$,
\begin{eqnarray*}
H(t) &=& H(t, R(t), U(t), \tilU(t))
      = \frac{1}{R(t)^2}
        \left| \left (\sigma(t, R(t) U(t)) - I \right)
               \left(U(t) - \tilU(t) \right)
        \right|,\\
G(t) &=& G(t, R(t), U(t))
      =  \frac{1}{R(t)^4}
        {\tr} \left( (\sigma(t, R(t) U(t)) - I)^2 \right).
\end{eqnarray*}
\end{lem}

{\it Proof.}
For a while let us use the notations,
\[
\tau = \tau(t, R(t), U(t)) = \sigma(t, R(t)U(t)) - U(t) \otimes U(t),
\]
\[
\tiltau = \tiltau (\tilU (t)) = I - \tilU (t) \otimes \tilU (t).
\]
Then the S.D.E.s governing $U(t)$
and $\tilU(t)$ can be written as
\[
dU = \frac{\tau}{R} dB - \frac{{\tr} (\tau^2)}{2R^2} U dt
\quad \mbox{and} \quad
d \tilU
= \frac{\tiltau}{R} d \tilB - \frac{{\tr} (\tiltau^2)}{2R^2} \tilU dt.
\]
Recall that the Brownian motions $B(t)$ and $\tilB (t)$
are adjusted so that
\[
dB = \omega(t) d \tilB (t)
\quad \mbox{and} \quad
U(t) = \omega(t) \tilU(t),
\]
where $\omega(t) = (\omega^{ij} (t))$ is the unitary matrix
which is given by
\[
\omega(t) = U(t) \otimes \tilU(t) - \tilU(t) \otimes U(t)
 + \langle U(t), \tilU(t) \rangle I.
\]
Consequently,
\begin{eqnarray*}
d \langle U, \tilU \rangle
&=& \langle \tilU, \frac{\tau}{R} dB \rangle
 - \frac{{\tr} ( \tau^2)}{2R^2} \langle \tilU, U \rangle dt
 + \langle U, \frac{\tiltau}{R} d \tilB \rangle\\
 && - \frac{{\tr} ( \tiltau^2)}{2R^2} \langle U, \tilU \rangle dt
 + \langle \frac{\tau}{R} dB, \frac{\tiltau}{R} d \tilB \rangle\\
&=& \frac{1}{R} \langle^t \omega \,^t \tau \tilU + \tiltau U, d\tilB \rangle
 - \frac{1}{2} R^2 G dt,
\end{eqnarray*}
where we put
\begin{equation}\label{GR2}
G = \frac{1}{R^2} \left\{
 \left( {\tr} (\tau^2) + {\tr}(\tiltau^2) \right)
 \langle U, \tilU \rangle
  - 2 {\tr} (\tiltau \tau \omega) \right\}.
\end{equation}
Soon we shall show that this quantity $G$ coincides with
the function $G(t)$ in the statement of the lemma.

Now let us compute the diffusion coefficients:
\begin{eqnarray*}
^t \omega \,^t\tau \tilU + \tiltau U
&=& ^t \omega(\sigma - I) \tilU + \,^t\omega \tilU
 - \,^t \omega \langle U, \tilU \rangle U + U - \langle U, \tilU \rangle \tilU\\
&=& ^t \omega(\sigma - I) \tilU + \,^t \omega \tilU + \omega \tilU
 - 2 \langle U, \tilU \rangle \tilU
= \, ^t \omega (\sigma - I) \tilU
\end{eqnarray*}
because $^t \omega + \omega = 2 \langle U, \tilU \rangle I$.
Hence we can define a 1-dimensional Brownian motion $W_1(t)$
through the relation
\[
\langle ^t\omega \,^t\tau \tilU + \tau \tilU, d\tilB \rangle
= \langle (\sigma - I) \tilU, dB \rangle
= R(t)^2 H dW_1,
\]
where we put
\[
H = \frac{1}{R(t)^2} | (\sigma - I) \tilU (t) |
  =  \frac{1}{R(t)^2} | (\sigma - I) (U(t) - \tilU (t)) |.
\]
Note that $dW_0 dW_0 = 0$.
In fact, by virtue of $(\sigma - I) U = 0$,
we have
\[
\langle (\sigma - I) \tilU, dB \rangle \langle U, dB \rangle
= \langle (\sigma - I) \tilU, U \rangle dt
= \langle \tilU, (\sigma - I) U \rangle dt
= 0.
\]
Consequently, the Brownian motions $W_0$ and $W_1$ are independent.

Lastly we check $G$ in (\ref{GR2})
coincides with $G(t)$ in the statement of the lemma.
Note that
\begin{eqnarray*}
&&\tr (\tiltau \tau \omega)
= \tr \left\{ \tiltau \tau \left( U \otimes \tilU - \tilU \otimes U
 + \langle U, \tilU \rangle I \right) \right\}\\
&=& \langle \tiltau \tau U, \tilU \rangle
 - \langle \tiltau \tau \tilU, U \rangle
 + \langle U, \tilU \rangle \tr(\tiltau \tau)\\
&=& \left( \langle \tau U, \tilU \rangle
         - \langle \tilU \otimes \tilU \tau U, \tilU \rangle
    \right)
  - \left( \langle \tau \tilU, U \rangle
    - \langle \tilU \otimes \tilU \tau U, U \rangle
    \right)\\
    &&+ \, \langle U, \tilU \rangle
     \tr \left( \tau (I - \tilU \otimes \tilU) \right)\\
&=& 0
 - \left( \langle (\sigma - U \otimes U) \tilU, U \rangle
 - \langle \tilU, \tau \tilU \rangle \langle \tilU, U \rangle \right)
 + \, \langle U, \tilU \rangle \left( \tr (\tau)
  - \langle \tau \tilU, \tilU \rangle \right)\\
&=& \langle U, \tilU \rangle \tr (\tau)
 =  \langle U, \tilU \rangle (\tr (\sigma) - 1).
\end{eqnarray*}
Hence we obtain
\begin{eqnarray*}
R^2 G
&=& \left\{ (\tr (\sigma^2) - 1) + (d - 1) \right\} \langle U, \tilU \rangle
 - 2 \langle U, \tilU \rangle \left( \tr (\sigma) - 1 \right)\\
&=& \left\{ \tr (\sigma^2) - 2 \tr (\sigma) + d \right\}
    \langle U, \tilU \rangle
= \tr ((\sigma - I)^2) \langle U, \tilU \rangle
\end{eqnarray*}
as is desired.
The proof is completed.

\section{The asymptotic analysis of the spherical parts}

Now we are ready to prove that the spherical motion
$U(t)$ is asymptotic to $\tilU(t)$ ``on the exponent".

\begin{prop}\label{prop:estimateUUtil}
For $U(t)$ and $\tilU(t)$ defined above,
the conditional expectation
\[
E\left[
 \exp \left\{ \frac{c}{\sqrt{\delta}} \left| U(t) - \tilU (t)\right| \right\}
    \left| \, \max_{[0, T]} R(t) < \delta \right. \right]
\]
remains bounded as $\delta$ tends to $0$ for every constant $c$.
\end{prop}

{\it Proof.}
Keeping in mind that $|U - \tilU|^2 = 2(1 - \langle U, \tilU \rangle)$,
let us consider the process
\begin{equation}\label{Delta1}
\Delta (t) = ( 1 - \langle U(t), \tilU(t) \rangle )
 \exp \left\{ \int_0^t \frac{1}{2} R(s)^2 G(s) ds \right\}.
\end{equation}
It is immediate to see that it satisfies the stochastic
differential equation
\begin{eqnarray*}
d \Delta (t)
&=& \left\{ - d \langle U, \tilU \rangle 
 + (1 - \langle U, \tilU \rangle ) \frac{1}{2} R^2 G dt \right\}
 \exp \int_0^t \frac{1}{2} R^2 G ds\\
&=& - R(t) H(t) \left\{ \exp \int_0^t \frac{1}{2} R^2 G ds \right\}
 dW_1(t) + d \left( \exp \int_0^t \frac{1}{2} R^2 G ds \right).
\end{eqnarray*}
Now let us take a 1-dimensional Brownian motions $W_2$,
which is independent of the Bessel process $R(t)$, such that
\[
- \int_0^t R(s) H(s) \left\{ \exp \int_0^s \frac{1}{2} R(r)^2 G(r) dr \right\}
 dW_1(s) = W_2 (\nu(t))
\]
where for simplicity of the notation we set
\[
\nu (t) = \int_0^t R(s)^2 H(s)^2 \left\{ \exp \int_0^s R(r)^2 G(r) dr
 \right\} ds.
\]
Therefore we can rewrite $\Delta(t)$ as follows:
\begin{equation}\label{Delta2}
\Delta(t) = W_2(\nu (t))
    + \left(\exp \int_0^t \frac{1}{2} R(s)^2 G(s) ds - 1\right).
\end{equation}
Notice that
\[
\nu (t) \leq \int_0^t R(s)^2 G(s) \left\{ \exp \int_0^s R(r)^2 G(r) dr
 \right\} ds
= \exp \int_0^t R(s)^2 G(s) ds - 1,
\]
since we have
\[
H(t)^2 = \frac{1}{R(t)^4} \left| (\sigma - I) \tilU (t)\right|^2
\leq
\frac{\tr ((\sigma - I)^2)}{R(t)^4} = G(t).
\]
Also notice that $G(t) < M$ for
a constant $M$ because $\sigma(t, x) - I = O(|x|^2)$.

Instead of the conditional expectation in the statement,
let us estimate the following conditional probability
\[
P_\delta (\lambda) = 
P \left[
 \sup_{[0, T]} \frac{1}{\sqrt{\delta}} \left|U(t) - \tilU(t)\right|
 > \lambda \left| \,
 \max_{[0, T]} R(t) < \delta \right. \right]
\]
for $\lambda >0$.
Since $| U - \tilU | / \sqrt{\delta} > \lambda$ is equivalent to
$2 ( 1 - \langle U, \tilU \rangle) / \delta > \lambda^2$,
by the definition (\ref{Delta1}) of $\Delta(t)$ we have
\begin{eqnarray*}
P_\delta (\lambda) &=&
P\left[ \sup_{[0, T]} \frac{2}{\delta}
      \, \Delta(t) \, e^{-\int_0^t (1/2) R^2(s) G(s) ds} > \lambda^2 \left| \,
    \max_{[0, T]} R(t) < \delta \right. \right]\\
&\leq&
P\left[ \sup_{[0, T]} \Delta(t) > \frac{\lambda^2 \delta}{2} \left| \,
    \max_{[0, T]} R(t) < \delta \right. \right].
\end{eqnarray*}
Substituting (\ref{Delta2}) for $\Delta (t)$,
\begin{eqnarray*}
P_\delta (\lambda)
&\leq&
 P\left[ \sup_{[0, T]} \left\{ W_2(\nu(t))
    + \left( e^{\int_0^t (1/2) R^2(s) G(s) ds} - 1 \right) \right\}
    > \frac{\delta \lambda^2}{2} \left| \,
    \max_{[0, T]} R(t) < \delta \right. \right]\\
&\leq&
 P\left[ \sup_{[0, T]} W_2(\nu(t))
    + \left( e^{(1/2) \delta^2 M T} - 1 \right)
    > \frac{\delta \lambda^2}{2} \left| \,
    \max_{[0, T]} R(t) < \delta \right. \right]\\
&\leq&
 P\left[ \sup_{[0, T]} W_2\left( e^{\int_0^t R(s)^2 G(s) ds} - 1 \right)
 > \frac{\delta \lambda^2}{2}
    - \left( e^{(1/2) \delta^2 M T} - 1 \right)
    \left| \, \max_{[0, T]} R(t) < \delta \right. \right]\\
&\leq&
 P\left[ \sup_{[0, \exp(\delta^2MT)-1]} W_2 (t)
 > \frac{\delta \lambda^2}{2}
    - \left( e^{(1/2) \delta^2 M T} - 1 \right)
     \left| \,
     \max_{[0, T]} R(t) < \delta \right. \right]\\
&=&
2P\left[ W_2(e^{\delta^2 MT} - 1)
  > \frac{\delta \lambda^2}{2}
    - \left( e^{(1/2) \delta^2 M T} - 1 \right)
    \right].
\end{eqnarray*}
We used the Andre reflection principle
and the independence of $W_2$ and $R$ in the last equality.
Therefore we have the final estimate
\[
P_\delta (\lambda)
\leq
2P\left[ W_2 (1) > \frac{\delta \lambda^2}{2}
    \left( e^{\delta^2 MT} - 1 \right)^{-1/2}
 - \left( e^{(1/2)\delta^2 M T} - 1 \right)
   \left( e^{\delta^2 M T} - 1 \right)^{-1/2} \right]
\]
by the scale invariance of the Brownian motion.
Since we can choose a constant $K$
independently on $\delta$
such that
\[
P_\delta (\lambda) \leq 2P \left[ W_2 (1) > K \lambda^2 \right]
\]
for large enough $\lambda$, the conditional expectation
\[
E \left[
 \exp \left\{ \frac{c}{\sqrt{\delta}} \left|U(t) - \tilU(t)\right| \right\} \left| \,
  \max_{[0, T]} R(t) < \delta \right. \right]
\]
is bounded in $\delta > 0$ for every constant $c$.
The proof is completed.

\section{The last step to the main theorem}

Now we can finish the proof of the main theorem.

By Proposition \ref{prop:preMain}
it suffices for the proof of Theorem \ref{thm:Main}
to show that
\[
E\left[ \exp \{ M(t) + L(t) \} \, \big| \, |Y|_T < \delta \right]
\to 1
\quad \mbox{as $\delta \to 0$}.
\]
In the following, we will show that its order is $\exp O(\sqrt{\delta} T)$.

Recall the following notations:
\[
M_p(t) = p \int_0^t \alpha_{ij} (s, Y(s)) dA^{ij} (s)
 - \frac{p^2}{2} \int_0^t \alpha_{ij} \alpha_{kl}
 d \langle A^{ij}, A^{kl} \rangle (s),
\quad M(t) = M_1(t),
\]
\[
L(t) = \int_0^t \beta(s, U(s)) ds
\quad \mbox{with} \quad
U(t) = \frac{Y(t)}{|Y(t)|},
\]
and 
\[
\beta(t, u) = \frac{d}{12} \ijdsum
 R_{ij} ( \gamma(t)) \left \{ u^i u^j - \frac{\delta^{ij}}{d} \right\}.
 \]
Set
\[
\tilde{L} (t) = \int_0^t \beta(s, \tilU(s)) ds
\]
where $\tilU (t) = \tilB (t) / |B(t)|$ is the spherical
Brownian motion introduced before.
Since $\beta (t, u)$ is Lipschitz with respect to $u$
uniformly in $t$ on $[0, T]$,
it follows from Proposition \ref{prop:estimateUUtil}
that there is a constant $C_1$ such that
\begin{equation}\label{L-Ltil}
E \left[
\exp \left. \left\{  \frac{1}{\sqrt{\delta}} (L(t) - \tilde{L}(t) \right\}
 \right| |Y|_T < \delta  \right] < C_1.
\end{equation}
We have also
that there exists a constant $C_2$ such that
\begin{equation}\label{Ltil}
E \left[ \left. \exp \left\{ \frac{1}{\sqrt{\delta}} \tilde{L}(T) \right\} \right|
 |Y|_T < \delta \right] < C_2.
\end{equation}

Now we apply H\"older inequality to our target expectation after
a little modification. First of all, note that for any positive
number $p$, we have
\begin{eqnarray}\label{M1Mp}
&&E \left[ \exp (M(t) + L(t)) \, \big| \, |Y|_T < \delta \right] \nonumber\\
&=& \exp (O(\delta T)) E\left[ \left. \exp \left\{
 \frac{1}{p} M_p (t) + L(t) \right\} \, \right| \, |Y|_T < \delta \right].
\end{eqnarray}

Now take small $\delta$ and set $p = (1 - 2\sqrt{\delta})^{-1}$.
Then, $1/p + \sqrt{\delta} + \sqrt{\delta} = 1$ and so
\begin{eqnarray}\label{Mp}
&&E \left[ \left. \exp \left\{ \frac{1}{p} M_p (t) + L(t)\right\} \right|
 |Y|_T < \delta \right] \nonumber\\
&\leq&
E \left[ \exp M_p(t) \, \big| \, |Y|_T < \delta \right]^{1/p}
\cdot E \left[ \left. \exp \left\{
    \frac{1}{\sqrt{\delta}} (L(t) - \tilde{L}(t))\right\}
\right| |Y|_T < \delta \right]^{\sqrt{\delta}} \nonumber\\
 &&\times \, E \left[ \left. \exp \left\{ \frac{1}{\sqrt{\delta}} \tilde{L}(t) \right\}
 \right| |Y|_T < \delta \right]^{\sqrt{\delta}}
\leq C_1^{\sqrt{\delta}} C_2^{\sqrt{\delta}} = \exp O(\sqrt{\delta} T).
\end{eqnarray}
Here we used the fact that $\exp M_p(t)$ is a martingale
under the conditional probability
$P[\, \cdot \, | \, |Y|_T < \delta]$.
From (\ref{M1Mp}) and (\ref{Mp}), we obtain the upper estimate
\[
E \left[ \exp \left\{ M(t) + L(t) \right\} \, \big| \, |Y|_T < \delta \right]
 < \exp O(\sqrt{\delta} T).
 \]
The inverse inequality is immediately obtained by Jensen inequality.
In fact,
\[
E \left[ \exp \left\{ M(t) + L(t) \right\} \, \big| \, |Y|_T < \delta \right]
\geq \exp E \left[ M(t) + L(t) \, \big| \, |Y|_T < \delta \right].
\]
It is obvious that
\[
\exp E \left[ M(t) \, \big| \, |Y|_T < \delta \right] = \exp O(\delta T).
\]
On the other hand, it follows from (\ref{Ltil}) and (\ref{L-Ltil})
using Proposition \ref{prop:estimateUUtil}
that
\begin{eqnarray*}
&&\exp E \left[ L(t) \, \big| \, |Y|_T < \delta \right]\\
&=& \exp E \left[ \tilde{L}(T) \, \big| \, |Y|_T < \delta \right]
 \exp E \left[ L(t) - \tilde{L}(t) \, \big| \, |Y|_T < \delta \right]\\
&=& \exp O(\sqrt{\delta} T) \exp O(\sqrt{\delta} T) = \exp O(\sqrt{\delta} T).
\end{eqnarray*}
Consequently,
\[
E \left[ \exp \left\{ M(t) + L(t) \right\} \, \big| \, |Y|_T < \delta \right]
\geq \exp O(\sqrt{\delta} T).
\]
Hence the proof of the Theorem is completed.

%
%

\begin{flushleft}
Keisuke HARA\\
Mynd, Inc.\\
Shirokanedai Bld. 3F\\
3-2-10 Shirokanedai, Minato-ku,\\
Tokyo, 108-0071 JAPAN\\
Email: hara.keisuke@gmail.com
\end{flushleft}

\end{document}